\documentclass[11pt,leqno]{article}
\usepackage{amsmath,amssymb}

\begin{document}

\hfil{\huge{\textbf{Existence globale de solutions}}}

\hfil{\huge{\textbf{d'{\'e}nergie infinie}}}

\hfil{\huge{\textbf{de l'{\'e}quation de Navier-Stokes 2D}}}

\bigskip

\hfil{\textbf{Pierre Germain}}

\hfil{\textbf{Centre de Math{\'e}matiques Laurent Schwartz, Ecole Polytechnique}}

\hfil{\textbf{U.M.R. 7640 du C.N.R.S.}}

\hfil{\textbf{91128 Palaiseau Cedex}}
\\

\newtheorem{theo}{Th{\'e}or{\`e}me}
\newtheorem{lemm}{Lemme}[section]
\newtheorem{prop}[lemm]{Proposition}
\newtheorem{rema}[lemm]{Remarque}

\textbf{R{\'e}sum{\'e}}

\bigskip

Nous {\'e}tudions dans cet article les solutions des {\'e}quations de
Navier-Stokes en deux dimensions, avec donn{\'e}e initiale dans
$\partial BMO$. Pour $u_{|t=0}$ dans l'adh{\'e}rence de la classe de
Schwartz, nous obtenons l'existence et l'unicit{\'e} d'une solution globale, et une
estimation sur sa norme dans $\partial BMO$.

\bigskip

\textbf{Abstract}

\bigskip

We study in this article the solutions of the Navier-Stokes equations,
with an initial data in $\partial BMO$. For $u_{|t=0}$ in the closure
of the Schwartz class, we obtain the existence and uniqueness of a
global solution, and an estimate on its norm in $\partial BMO$.

\bigskip

\section{Introduction}

\subsection{Les {\'e}quations de Navier-Stokes}

Nous nous int{\'e}resserons dans cet article aux solutions globales du
probl{\`e}me de Cauchy associ{\'e} aux {\'e}quations de Navier Stokes, qui
d{\'e}crivent le mouvement d'un fluide visqueux dans l'espace tout
entier. Ce mouvement est
d{\'e}crit par les variables $u(x,t)$ et $p(x,t)$ qui donnent, en un
point $x$ de l'espace et au temps $t$, respectivement la vitesse et la
pression du fluide. Nous consid{\'e}rons le cas d'un fluide
incompressible (ce qui entra{\^\i}ne la condition $\operatorname{div} u =
0$), et de viscosit{\'e} $\nu = 1$. Les {\'e}quations de Navier-Stokes
prennent alors la forme suivante~:
\begin{equation*}
 (NS) \left\{ \begin{array}{l} 
\partial_{t}u - \Delta u + u \cdot \nabla u = -\nabla p \\
\operatorname{div} u = 0 \\
u_{|t=0} = u_{0} \, \, .
\end{array} \right.
\end{equation*}
Nous consid{\'e}rerons aussi bri{\`e}vement le cas o{\`u} une force
ext{\'e}rieure d{\'e}rivant d'un potentiel~$V$ est appliqu{\'e}e. Le
syst{\`e}me devient alors~:
\begin{equation*}
 (NSF) \left\{ \begin{array}{l} 
\partial_{t}u - \Delta u + u \cdot \nabla u = -\nabla p + \nabla V \\
\operatorname{div} u = 0 \\
u_{|t=0} = u_{0} \, \, .
\end{array} \right.
\end{equation*}
Sauf mention contraire, nous nous placerons toujours \textbf{en deux
  dimensions d'espace}. 

\subsection{Solutions de Leray}

La th{\'e}orie des {\'e}quations de Navier-Stokes a {\'e}t{\'e} initi{\'e}e par
Jean Leray~; il s'est int{\'e}ress{\'e} aux solutions d'{\'e}nergie finie de
$(NS)$ issues de $u_{0} \in L^2$.

L'espace $L^2$ est l'espace d'{\'e}nergie pour les conditions initiales~;
l'espace d'{\'e}nergie associ{\'e} aux solutions de l'{\'e}quation de
Navier-Stokes est d{\'e}fini par~:
$$ \mathcal{L} \overset{\mbox{d{\'e}f}}{=} L^{\infty}(\mathbb{R}^{+},L^{2}) \cap
L^{2}(\mathbb{R}^{+},\dot{H}^{1}) \, \, ,$$
o{\`u} $\dot{H}^1$ est l'espace de Sobolev homog{\`e}ne, soit l'ensemble
des fonctions $f$ pour lesquelles
$$
\int_{\mathbb{R}^n} |\xi|^2 |\hat{f}(\xi)|^2 \, d\xi \, < \, \infty
\, \, .
$$
La norme de $\mathcal{L}$ s'{\'e}crit :
$$ \| v \|_{\mathcal{L}} = \| v
\|_{L^{\infty}(\mathbb{R}^+,L^{2})} + \| v 
\|_{L^{2}(\mathbb{R}^+,\dot{H}^{1})}
. $$

Les solutions de $(NS)$ appartenant {\`a} cet espace sont appel{\'e}es
solutions de Leray. Le th{\'e}or{\`e}me suivant garantit leur existence et
leur unicit{\'e}.

\begin{theo}[J. Leray~\cite{bibleray}]
\label{theoleray}
Soit $v_{0} \in L^{2}(\mathbb{R}^{2})$. Il existe une unique solution
globale $v~\in~\mathcal{L}$ de $(NS)$ ayant
$v_{0}$ pour condition initiale.
\end{theo}

\subsection{Equation int{\'e}grale, {\'e}quation projet{\'e}e}
La formulation int{\'e}grale des {\'e}quations de Navier-Stokes $(NS)$ s'{\'e}crit~:
\begin{equation}
\label{NSint}
\begin{split}
u(t) &= e^{t\Delta}u_{0} - \int_{0}^{t}e^{(t-s)\Delta} \mathbb{P}
\nabla \cdot (u(s) \otimes u(s))ds \\
&= e^{t\Delta}u_{0} - B(u,u) \, \, , \\
\end{split}
\end{equation}
o{\`u} l'on note $\mathbb{P}$ le projecteur de Leray sur les champs
de vecteur {\`a} divergence nulle, et $B$ l'op{\'e}rateur bilin{\'e}aire
$$
B(u,v) \overset{\mbox{d{\'e}f}}{=} \int_{0}^{t}e^{(t-s)\Delta}
  \mathbb{P} \nabla \cdot (u(s) \otimes v(s))ds \, \, .
$$

\medskip

On peut montrer que le noyau de $\nabla^{k} e^{t\Delta} \mathbb{P}
\nabla \cdot$ s'{\'e}crit~: 
$$ t^{-\frac{3+k}{2}}G_{k}(\frac{|x|}{\sqrt{t}}) $$
avec $G_{k} \in L^{1}\cap L^{\infty}$. On notera pour plus de
simplicit{\'e} indiff{\'e}remment $G$ pour tous les~$G_{k}$.
De m{\^e}me, nous consid{\'e}rerons dans ce qui suit pour all{\'e}ger les
{\'e}critures que $B$ op{\`e}re sur des fonctions
r{\'e}elles. Ceci nous conduit aux notations suivantes~\cite{bibcannone}~:
\begin{gather*}
B(u,v) = \int_{0}^{t} \frac{1}{(t-s)^{3/2}}
  G \left( \frac{\cdot}{\sqrt{t-s}} \right) \ast (u(s) v(s))
  ds \\
\nabla B(u,v) = \int_{0}^{t} \frac{1}{(t-s)^{2}} 
  G \left( \frac{\cdot}{\sqrt{t-s}} \right) \ast (u(s) v(s)) ds
\end{gather*}

\bigskip

Suivant un abus de langage habituel, nous appellerons dans cet article
"solutions de $(NS)$" des \textbf{solutions de l'{\'e}quation
  int{\'e}grale}~(\ref{NSint}) (solutions "\textit{mild}"). Les
solutions de~(\ref{NSint}) ne sont pas a priori des solutions du
syst{\`e}me $(NS)$ pr{\'e}sent{\'e} en t{\^e}te de cet
article, voir~\cite{biblemarie} pour une discussion {\`a} ce sujet. Par
contre, pour ce qui est des solutions que nous consid{\'e}rerons,
l'{\'e}quation int{\'e}grale est {\'e}quivalente {\`a} l'{\'e}quation
projet{\'e}e :
$$
\left\{ \begin{array}{l} 
\partial_{t}u - \Delta u + \mathbb{P}(u \cdot \nabla u) = 0 \\
\operatorname{div} u = 0 \\
u_{|t=0} = u_{0}
\end{array} \right.
$$
Cette {\'e}quation sera consid{\'e}r{\'e}e au sens des distributions. Notons
que l'on peut (au moins formellement) revenir {\`a} $(NS)$ en calculant
la pression gr{\^a}ce {\`a} la formule
$$
-\Delta p = \operatorname{div}(u\cdot \nabla u) \; .
$$
\subsection{Espaces critiques pour les conditions initiales}

La d{\'e}finition des espaces critiques repose sur les propri{\'e}t{\'e}s
d'invariance par dilatation et translation des solutions de $(NS)$. Soyons plus
explicites~: soit $u$ une solution de $(NS)$ associ{\'e}e {\`a} la condition
initiale $u_{0}$. Alors, si $\lambda >0$ et $x_{0} \in
\mathbb{R}^{2}$, $\lambda u(\lambda(x-x_{0}),\lambda^2 t)$ sera
associ{\'e}e {\`a} $\lambda u_{0} (\lambda ( x-x_{0}))$. Un espace de Banach
$X\hookrightarrow \mathcal{S}'$ sera dit critique (sous-entendu~: pour
les conditions initiales de $(NS)$) si sa norme v{\'e}rifie~: 
\begin{equation}
\label{invechelle}
\forall \lambda>0 \, , \, \forall x_0 \in \mathbb{R}^2 \, \, \, \| u \|_X = \lambda
\| u (\lambda (\cdot - x_{0})) \|_X \, \,.
\end{equation}
On voit facilement que, en dimension deux d'espace, $L^2$, les espaces
de Besov $\dot{B}^{\frac{2}{p}-1}_{p,q}$ (avec $1\leq p,q\leq \infty$)
et  $\partial BMO$ sont des espaces critiques (se reporter {\`a} la
troisi{\`e}me partie pour une d{\'e}finition de ces espaces). 

\bigskip

Il est int{\'e}ressant de noter (\cite{bibdub}) que tout espace critique
$X$ v{\'e}rifie~:
$$
X\hookrightarrow \dot{B}_{\infty,\infty}^{-1} \, .
$$
D'autre part (\cite{bibdub}, \cite{bibkoch}), on a la cha{\^\i}{}ne d'inclusions suivante~:
\begin{equation}
\label{chaininclu}
L^{2} \hookrightarrow \dot{B}^{\frac{2}{p}-1}_{p,q}
\hookrightarrow \dot{B}^{\frac{2}{p}-1}_{p,\infty} 
\hookrightarrow \partial BMO \hookrightarrow \dot{B}^{-1}_{\infty,\infty}
\end{equation}
o{\`u} $p,q \in [2,\infty[$. 

\bigskip

En dimension d'espace quelconque, des th{\'e}or{\`e}mes garantissent
l'existence d'une solution locale $u^\star$ pour $(NS)$ si $u_0$ est
de norme grande dans un espace critique $X$~: voir~\cite{bibchemin}
pour~$\dot{B}^{\frac{2}{p}-1}_{p,q}$ et \cite{bibdub}, \cite{bibkoch}
pour $\partial BMO$.

En dimension 2 d'espace, nous avons vu que l'espace d'{\'e}nergie
$L^2$ est lui-m{\^e}me un espace critique. En utilisant cette
propri{\'e}t{\'e}, I. Gallagher et F. Planchon~\cite{bibgallag} ont pu
montrer que, pour une donn{\'e}e initiale grande dans
$\dot{B}^{\frac{2}{p}-1}_{p,q}$, la solution locale en temps $u^\star$
se prolonge {\`a}~$\mathbb{R}^+$, et obtenir ainsi une solution globale
en temps. Nous nous proposons d'appliquer la m{\^e}me m{\'e}thode {\`a}
$\partial BMO$.

\bigskip

Mais avant de nous int{\'e}resser {\`a} $(NS)$ avec $u_0 \in \partial BMO$, deux
remarques doivent {\^e}tre faites sur cet espace~:
\begin{itemize}
\item Un probl{\`e}me pos{\'e} par $\partial BMO$ est que
  $\partial BMO \neq \overline{\mathcal{S}}^{\partial BMO}$ (l'adh{\'e}rence de la
  classe de Schwartz dans $\partial BMO$)~; on retrouve ce probl{\`e}me pour
  les espaces de Besov d'indice~$q$ infini. Nous
  nous restreindrons donc souvent dans ce qui suit {\`a}
  $\overline{\mathcal{S}}^{\partial BMO}$.
\item Signalons d'autre part que $\partial BMO$ est, en un certain sens,
  l'espace optimal pour la r{\'e}solution de $(NS)$.
En effet, $\partial BMO$ est le plus gros espace critique
$X$, dans la cha{\^\i}{}ne d'inclusions pr{\'e}c{\'e}dente (qui se
g{\'e}n{\'e}ralise {\`a} $\mathbb{R}^n$), pour lequel on sait
d{\'e}montrer un r{\'e}sultat d'existence de solutions de $(NS)$ dans
$\mathbb{R}^n$ avec $u_{0} \in X$ petit. En fait, on peut formaliser et fonder
cette observation, et donner un sens pr{\'e}cis {\`a} l'id{\'e}e que~$\partial
  BMO$ est optimal~: voir~\cite{bibauscher}.
\end{itemize}

\subsection {Solutions de Koch et Tataru}

H. Koch et D. Tataru ont, dans~\cite{bibkoch}, obtenu l'existence de
solutions globales de $(NS)$, en toute dimension, pour une donn{\'e}e
initiale petite dans $\partial BMO$. Avant d'{\'e}noncer ce r{\'e}sultat plus
pr{\'e}cis{\'e}ment, il nous faut introduire les espaces $X_{T}$, $0<T\leq
\infty$~; $X_{T}$ 
est l'ensemble des fonctions $w$ pour lesquelles la norme suivante est
bien d{\'e}finie et finie~:
\begin{equation*}
\begin{split}
\| w \|_{X_{T}} & \overset{\mbox{d{\'e}f}}{=}
\sup_{0<t<T} \sqrt{t} \| w(t)
\|_{\infty} + \sup_{0<t<T} t \| \nabla w(t)
\|_{\infty} \\
& + \sup_{0<R<\sqrt{T}} \sup_{x\in \mathbb{R}^{2}}
\left( \int_{0}^{R^{2}} {\int\!\!\!\!\!\!-}_{B(x,R)} \mid w(y,t)
  \mid^{2} dy\, dt 
\right)^{1/2} \, \, ,
\end{split}
\end{equation*}
o{\`u} l'on note pour la moyenne de $f$ sur la boule $B(x,R)$~: $\displaystyle
{\int\!\!\!\!\!\!-}_{B(x,R)} f = \frac{1}{|B(x,R)|}\int_{B(x,R)}
f$. Nous verrons dans la troisi{\`e}me partie que le troisi{\`e}me terme
de cette norme est adapt{\'e} {\`a} des donn{\'e}es initiales dans $\partial
BMO$ puisqu'il est fini si $w=e^{t\Delta}w_0$ avec $w_0 \in \partial BMO$.

\bigskip

Les d{\'e}monstrations des r{\'e}sultats suivants peuvent {\^e}tre trouv{\'e}es
dans~\cite{bibkoch} et~\cite{bibdub}.

\begin{theo}[H. Koch, D. Tataru, \cite{bibkoch}]
\label{theokochtataru}
Soit $d \geq 2$. Il existe $\epsilon>0$ tel que, si $w_{0} \in
\partial BMO(\mathbb{R}^d)$ est de divergence nulle et
$\| w_{0} \|_{\partial BMO} < \epsilon$, il existe $w$ une
solution globale de~$(NS)$ v{\'e}rifiant l'in{\'e}galit{\'e} suivante~:
\begin{equation} 
\label{estimtataru} \| w \|_{X_{\infty}} \leq
C \| w_0 \|_{\partial BMO} \, \, ,
\end{equation}
pour une constante C. De plus, si $w_{0} \in \overline{\mathcal{S}}^{\partial BMO}$, alors 
$$\|
w \|_{X_{t}} \underset{t \rightarrow 0}{\longrightarrow} 0 \; .$$
\end{theo}

\textsc{Id{\'e}e de la preuve~:}
On consid{\`e}re l'{\'e}quation int{\'e}grale~(\ref{NSint}), {\`a} laquelle on
applique le th{\'e}or{\`e}me de point fixe~\ref{pointfixe}, en se pla{\c c}ant dans
l'espace $X_{\infty}$,
d{\'e}fini ci-dessus. On s'assure donc des points suivants :
\begin{itemize}
\item $\| e^{t\Delta}w_{0} \|_{X_{\infty}} \leq \|
  w_{0} \|_{\partial BMO}$
\item $B~:X_{\infty} \times X_{\infty} \rightarrow X_{\infty}$ est
  bicontinue~; ceci est garanti par la proposition suivante~:
\end{itemize}

\begin{prop}[\cite{bibkoch}]
\label{propBbicont}
Il existe $\eta > 0$ tel que, $\forall T \in ]0,\infty]$, 
\begin{equation}
\label{Bbicont}
\| B(u,v) \|_{X_{T}} \leq \eta \| u
\|_{X_{T}} \| v \|_{X_{T}} \, .
\end{equation}
On dispose aussi de l'estimation, si $0<t<T$, et pour une constante
$C>0$~:
\begin{equation}
\label{BdBMO}
\| B(u,v) (t) \|_{\partial BMO} \leq C \| u \|_{X_T} \| v
\|_{X_T} \, .
\end{equation}
\end{prop}

On peut donc appliquer le th{\'e}or{\`e}me~\ref{pointfixe}, et le
th{\'e}or{\`e}me~\ref{theokochtataru} en d{\'e}coule. $\blacksquare$

Le th{\'e}or{\`e}me pr{\'e}c{\'e}dent s'adapte facilement pour fournir
l'existence de solutions locales {\`a} donn{\'e}es grandes~; cependant on doit
se restreindre
{\`a} des donn{\'e}es dans $\overline{\mathcal{S}}^{\partial BMO}$, pour pouvoir
garantir que $\| e^{t\Delta} w_{0} \|_{X_{T}}
\underset{T\rightarrow 0}{\longrightarrow} 0$ et appliquer le th{\'e}or{\`e}me de point
fixe~\ref{pointfixe} {\'e}nonc{\'e} en annexe. On obtient alors~:

\begin{theo}
\label{KTdonneesgrandes}
Soit $u_{0} \in \overline{\mathcal{S}}^{\partial BMO}(\mathbb{R}^d)$ de
divergence nulle. Il existe $T^{\star}>0$ et $u
\in X_{T^{\star}}$ une solution de $(NS)$ sur $[0,T^{\star}]$ issue de
$u_{0}$. De plus, $\| u \|_{X_{t}}
  \underset{t \rightarrow 0}{\longrightarrow} 0$.
\end{theo}

Une question naturelle est maintenant celle de l'unicit{\'e} des
solutions de Koch et Tataru. Soit $\mathcal{E}_{T}$ l'espace
\begin{equation}
\label{unptiesp}
\mathcal{E}_{T} \overset{\mbox{d{\'e}f}}{=} \left\{ f \in
L^{\infty}_{\operatorname{loc}}(]0,T],L^{\infty})\mbox{ , }
  \underset{\delta \rightarrow 0}{\operatorname{lim}} \| f
  \|_{X_{\delta}}<\frac{1}{2\eta}
\right\} \, .
\end{equation}

\begin{prop}[\cite{bibdub}]
\label{propunicite}
Deux solutions de $(NS)$ appartenant {\`a} l'espace $\mathcal{E}_{T}$ pour
un $T>0$ et issues d'un m{\^e}me $u_{0} \in \partial BMO$ sont {\'e}gales.
\end{prop}

\section{ Etude de $(NS)$ et $(NSF)$ pour $u_{0} \in
  \overline{\mathcal{S}}^{\partial BMO}$}

\subsection{R{\'e}sultat principal}

Nous nous proposons dans le pr{\'e}sent article d'{\'e}tendre {\`a} des
donn{\'e}es grandes dans $\partial BMO$, et en deux dimensions, le r{\'e}sultat
d'existence globale (en temps) obtenu par H. Koch et D. Tataru pour
des donn{\'e}es petites dans
$\partial BMO$.

Il nous faut ici faire mention de~\cite{bibgiga}. Les auteurs prouvent
dans cet article par des m{\'e}thodes d'analyse harmonique l'existence, pour $u_{0} \in
L^{\infty}(\mathbb{R}^2)$, d'une solution globale de $(NS)$ dans
$L^{\infty}_{\operatorname{loc}}([0,\infty[,L^{\infty})$. Or le
th{\'e}or{\`e}me~\ref{KTdonneesgrandes} donne pour une condition initiale
dans~$\overline{\mathcal{S}}^{\partial BMO}$ l'existence locale d'une
solution de $(NS)$ de norme $L^{\infty}$ finie pour $t>0$. Ces deux
{\'e}l{\'e}ments permettent de d{\'e}duire, pour $u_{0} \in
  \overline{\mathcal{S}}^{\partial BMO}$, l'existence d'une solution
  globale de~$(NS)$.

Cependant, on a oubli{\'e} ce faisant le cadre $\partial BMO$. La
m{\'e}thode que nous pr{\'e}sentons ici, outre qu'elle prouve l'existence
d'une solution globale, permet aussi de montrer que cette solution est
dans $L^{\infty}_{\operatorname{loc}}([0,\infty[, \partial BMO)$, et donne
une estimation sur sa norme dans $\partial BMO$.

\bigskip

Avant d'{\'e}noncer le th{\'e}or{\`e}me, il nous faut d{\'e}finir l'espace
dans lequel l'unicit{\'e} sera prouv{\'e}e~; nous aurons besoin de la norme de
type Carleson suivante (cette norme intervient d{\'e}j{\`a} dans la
d{\'e}finition de $X_T$)~:
\begin{equation*}
\| u \|_{\mathcal{C},T} = \sup_{0<R<\sqrt{T}} \sup_{x\in \mathbb{R}^{2}}
\left( \int_{0}^{R^{2}} {\int\!\!\!\!\!\!-}_{B(x,R)} | u(y,t)
  |^{2} dy\, dt \right)^{1/2} \, \, .
\end{equation*}
Nous pouvons maintenant d{\'e}finir l'espace d'unicit{\'e}
\begin{equation}
\label{defdeE}
u \in \mathcal{E} \overset{\mbox{d{\'e}f}}{\Longleftrightarrow} 
\left\{ \begin{array}{l}
(i) \, \, \, \| u \|_{\mathcal{C},T} < \infty \, \, \forall T \\
(ii) \, \, \, \forall
T \in \mathbb{R}^+ , \exists \epsilon > 0 \mbox{ tel que } u \in
L^\infty_{loc} (]T,T+\epsilon[,L^\infty) \\ 
\; \; \; \; \; \; \; \; \; \; \; \; \; \; \; \; \; \; \; \; \; \; \;
\; \; \; \; \; \; \; \; \; \; \; \; \; \; \; \; \; \;

\; \mbox{ et } \underset{\delta
  \rightarrow 0}{\operatorname{lim}} \| u(T+\cdot)
\|_{X_\delta} < \frac{1}{2\eta}
\end{array} \right.
\; 
\end{equation}
o{\`u} $\eta$ est la constante strictement positive intervenant dans la
proposition~\ref{propBbicont} et dans la d{\'e}finition de
$\mathcal{E}_T$. La condition $(i)$ va nous
permettre de garantir, si en
outre $u$ est solution de $(NS)$, que $u$ est faiblement
continue. Quant {\`a} $(ii)$, elle signifie que pour tout $\tau
\in \mathbb{R}^+$, $u(\tau+.)$
appartient {\`a} $\mathcal{E}_\epsilon$ (d{\'e}fini en~(\ref{unptiesp})) pour un certain $\epsilon > 0$.

\begin{rema}
Notons que la classe d'unicit{\'e} qui vient d'{\^e}tre d{\'e}finie englobe
des r{\'e}sultats d{\'e}j{\`a} connus sur $(NS)$ en deux dimensions, avec
$u_{0}$ dans $\partial BMO$ ou un sous-espace de~$\partial BMO$.
\begin{itemize}
\item
On v{\'e}rifie ais{\'e}ment que notre classe d'unicit{\'e} comprend les
solutions $L^{\infty}_{loc}(\mathbb{R}^+,L^{\infty})$ construites 
dans~\cite{bibgiga}, car ces derni{\`e}res sont localement
born{\'e}es dans $L^\infty$.
\item
Elle comprend aussi les solutions de Koch et
Tataru (th{\'e}or{\`e}mes~\ref{theokochtataru} et~\ref{KTdonneesgrandes}) :
si $u$ est l'une de ces solutions, elle
v{\'e}rifie clairement $(i)$. Elle v{\'e}rifie d'autre part $(ii)$ en~$\tau =
0$~; en $\tau > 0$ aussi car $u$ est born{\'e}e dans $L^\infty$ sur
$[\delta, T[$ pour tout $\delta >0$, si l'on note $T$ le temps
d'existence de la solution.
\item Elle comprend enfin les solutions uniques dans
  $\mathcal{C}(\mathbb{R}^+,\dot{B}_{p,q}^{-1+\frac{2}{p}})$, avec
  $q<\infty$ construites dans~\cite{bibgallag}. 

[En effet, soit $u$ une solution dans
  $\mathcal{C}(\mathbb{R}^+,\dot{B}_{p,q}^{-1+\frac{2}{p}})$. L'hypoth{\`e}se $q<\infty$ implique que $u$ est {\`a}
  valeurs dans $\overline{\mathcal{S}}^{\partial BMO}$. On se place
maintenant au temps $\tau \geq 0$, puis on r{\'e}soud
  l'{\'e}quation par point fixe, pour $T$ assez petit, dans $X_{T} \cap
  \tilde{L}^{r}([0,T],\dot{B}^{-1+\frac{2}{p} + \frac{2}{r}}_{p,q})$
  (voir~\cite{bibgallag2}), et on conclut par unicit{\'e} de la solution
  d'un probl{\`e}me de point fixe que $u(\tau+\cdot)$ appartient pour~$T$ assez petit {\`a}
  $X_{T}$. Ceci implique que $(ii)$ est v{\'e}rifi{\'e}e, et que $(i)$
  l'est aussi pour $T$ proche de $0$.

D'autre part, $u$ s'{\'e}crit (se reporter {\`a}~\cite{bibgallag}) $u =
  v+w$, avec $w$ born{\'e}e dans $L^\infty$ et $v$ localement d'{\'e}nergie
  finie pour $t>0$, donc $(i)$ est v{\'e}rifi{\'e}e pour tout $T$. ]
\end{itemize}
\end{rema}

Il est temps d'{\'e}noncer notre th{\'e}or{\`e}me principal.

\begin{theo}
\label{theoprincipal}
Soit $u_{0} \in \overline{\mathcal{S}}^{\partial BMO}(\mathbb{R}^2)$ de
  divergence nulle. Il existe une solution globale de~$(NS)$, $u$, {\`a} valeurs
  dans $\partial BMO$ issue de
  $u_{0}$. Cette solution v{\'e}rifie pour tout $\delta >0$ et pour une constante $C(\delta)$ d{\'e}pendant de 
  $u_{0}$ et de $\delta$~:
$$
\forall t>0 \, \, \, \, \| u(t) \|_{\partial BMO} \leq C(\delta) (1+t^\delta) \, .
$$
De plus, cette solution prolonge la solution de Koch et Tataru
(th{\'e}or{\`e}me~\ref{KTdonneesgrandes}), d{\'e}finie seulement localement, et est
unique dans la classe $\mathcal{E}$.
\end{theo}

La preuve de ce th{\'e}or{\`e}me est l'objet des sections~\ref{preuvedutheo1} et~\ref{preuvedutheo2}.

\begin{rema}
En fait, nous prouvons dans ce qui suit que $u$ peut s'{\'e}crire $u =
v+w$, avec $\underset{t>0}{\operatorname{sup}} \sqrt{t} \| w(t)
\|_{\infty}$ born{\'e}, $\sqrt{t} \| v(t) \|_{\dot{H}^1}$ born{\'e} pr{\`e}s
de 0, et $v \in L^{2}_{loc} \dot{H}^{1}$ en dehors de 0. Les
inclusions classiques $L^\infty \hookrightarrow BMO$ et $\dot{H}^{1}
\hookrightarrow BMO$ impliquent que $u$ appartient {\`a} l'espace
$$
L^\infty_{\operatorname{loc}} (\mathbb{R}^+,\partial BMO) \cap
L^{2,\infty}_{\operatorname{loc}} (\mathbb{R}^+,BMO) \, \, .
$$
Cet espace pr{\'e}sente une analogie claire avec
$$
L^\infty (\mathbb{R}^+,L^2) \cap
L^{2} (\mathbb{R}^+,\dot{H}^1) \, \, ,
$$
qui est l'espace d'{\'e}nergie $\mathcal{L}$ d{\'e}fini plus haut.
\end{rema}

\subsection{Construction de la solution}

\label{preuvedutheo1}
Cette section est consacr{\'e}e {\`a} la construction de la solution qui
fait l'objet du th{\'e}or{\`e}me~\ref{theoprincipal}~; il s'agit donc de
la preuve de la partie "existence" du th{\'e}or{\`e}me~\ref{theoprincipal}.

\medskip

Soit $u_{0} \in
\overline{\mathcal{S}}^{\partial BMO}$. Nous proc{\'e}dons de la mani{\`e}re
suivante~:
\begin{enumerate}
\item $u_{0}$ s'{\'e}crit comme somme d'une partie grande  et r{\'e}guli{\`e}re
$v_{0}$, et d'une partie petite et moins r{\'e}guli{\`e}re $w_{0}$. Le
th{\'e}or{\`e}me de Koch et Tataru donne l'existence d'une solution~$w$ de
$(NS)$ issue de $w_{0}$~; on se ram{\`e}ne ainsi {\`a} une
{\'e}quation (faisant intervenir $w$) d'inconnue $v$, et de donn{\'e}e initiale $v_{0}$. 
\item L'utilisation
d'un th{\'e}or{\`e}me de point fixe permet alors d'obtenir une solution $v$
appartenant {\`a} $X_T$ et locale en temps.
\item Par un argument de propagation de r{\'e}gularit{\'e}, on montre
 qu'il existe un temps strictement positif pour lequel $v \in L^2$.
\item Une estimation d'{\'e}nergie a priori permet de rendre $v$ globale
 en temps.
\end{enumerate}
C'est la m{\'e}thode utilis{\'e}e par I. Gallagher et F. Planchon
dans~\cite{bibgallag} pour des donn{\'e}es initiales dans des espaces de
Besov. Cette id{\'e}e a {\'e}t{\'e} utilis{\'e}e pour la premi{\`e}re fois dans
le cadre des {\'e}quations de Navier-Stokes par
Calder{\'o}n~\cite{bibcalderon}.

\bigskip

Comme nous l'avons vu, la premi{\`e}re {\'e}tape de la d{\'e}monstration consiste
{\`a} d{\'e}couper la donn{\'e}e initiale. Soit donc $u_{0} \in
\overline{\mathcal{S}}^{\partial BMO}$ de divergence nulle. Par
d{\'e}finition de cet espace, il 
existe $v_{0}$ et $w_{0}$, tous deux de divergence nulle, et tels que~:
\begin{itemize}
\item $u_0 = v_0 + w_0$
\item $v_{0} \in \mathcal{S}$
\item $w_{0} \in \partial BMO$, et $\| w_{0} \|_{\partial BMO}
  \leq \epsilon$ ($\epsilon$ pourra {\^e}tre fix{\'e} ult{\'e}rieurement. On
  le prendra assez petit pour que tous les
  th{\'e}or{\`e}mes dont nous aurons besoin s'appliquent).
\end{itemize}

Par le th{\'e}or{\`e}me~\ref{theokochtataru}, il existe $w$ une solution de
$(NS)$ issue de $w_{0}$. De plus, $w$ v{\'e}rifie les
estimations~(\ref{estimtataru}).
Si $u$ v{\'e}rifie l'{\'e}quation~(\ref{NSint}), on voit ais{\'e}ment que
$v$ est solution de~:
\begin{equation}
\label{equationenv}
v(t) = e^{t\Delta}v_{0} - B(v,w)(t) - B(w,v)(t) - B(v,v)(t)
\end{equation}

Nous avons maintenant besoin de la d{\'e}finition suivante~: soit $Y_{T}$
l'espace d{\'e}fini par~:
$$
\| f \|_{Y_{T}} \overset{\mbox{d{\'e}f}}{=} \| f
\|_{L^{\infty}([0,T],L^{2})} + \sup_{0<t<T} \sqrt{t} \|
  \nabla f(t) \|_{L^{2}} \, \, .
$$

La proposition suivante va nous donner l'existence d'une
solution locale dans $X_{T}\cap Y_{T}$.

\bigskip

\begin{prop}
Il existe $T>0$ tel que~(\ref{equationenv}) admette une solution $v$
sur $[0,T]$ appartenant {\`a} $X_{T} \cap Y_{T}$.
\end{prop}

\textsc{Preuve de la proposition~:} Nous commen{\c c}ons
par construire une solution dans $X_{T}$~; nous montrerons ensuite que
cette solution appartient aussi {\`a} $Y_{T}$. Il s'agit dans un premier
temps de pouvoir appliquer le th{\'e}or{\`e}me de point fixe~\ref{pointfixe}
(voir l'annexe) {\`a} la
r{\'e}solution dans~$X_{T}$ de~(\ref{equationenv}), c'est {\`a} dire de
s'assurer qu'il existe un $T$ tel que le th{\'e}or{\`e}me s'applique. Il
nous faut v{\'e}rifier les points suivants~:
\begin{itemize}
\item $\| e^{t\Delta}v_{0} \|_{X_{T}} \underset{T
    \rightarrow 0}{\longrightarrow} 0$~; c'est le cas pour $v_{0} \in
    \mathcal{S}$.
\item $B~:X_{T} \times X_{T} \rightarrow X_{T}$ est bicontinue~; c'est
  vrai d'apr{\`e}s~(\ref{Bbicont}).
\item $\| B(w,\cdot) \|_{\mathcal{L}(X_{T})}$ et $\|
  B(\cdot,w) \|_{\mathcal{L}(X_{T})}$ sont finies et leur somme est
  strictement inf{\'e}rieure {\`a} 1 ; c'est
  le cas si $\| w \|_{X_{T}}$ est assez petit, du fait
  de~(\ref{Bbicont}).
\end{itemize}

\bigskip

Nous obtenons ainsi $T>0$ et $v$ solution de~(\ref{equationenv}) sur
$[0,T]$. Nous allons maintenant montrer que $v$ appartient aussi {\`a}
$Y_{T}$, en appliquant le lemme de propagation de
r{\'e}gularit{\'e}~\ref{propagation} (avec, pour
reprendre les notations du lemme, $X=X_{T}$ et $Y=Y_{T}$) {\`a}
l'{\'e}quation~(\ref{equationenv}). Il nous faut nous assurer des points suivants~:
\begin{itemize}
\item $e^{t\Delta}v_{0} \in Y_{T}$~; ceci est v{\'e}rifi{\'e} d{\`e}s que
  $v_{0} \in \mathcal{S}$. C'est ici que cette hypoth{\`e}se intervient.
\item $B~:X_{T}\times Y_{T}\rightarrow Y_{T}$ et $B:Y_{T}\times
  X_{T}\rightarrow Y_{T}$ sont bicontinues~; c'est l'objet de la
  proposition suivante, dont nous renvoyons la
  d{\'e}monstration {\`a} la fin de cette section.
\end{itemize}
\begin{prop}
\label{continuitedutruc}
L'application $B$ est bicontinue de $X_{T} \times Y_{T} \rightarrow
Y_{T}$.
\end{prop}
\begin{itemize}
\item $B(w,\cdot)~:Y_{T} \rightarrow Y_{T}$ et $B(\cdot,w)~:Y_{T}
  \rightarrow Y_{T}$ sont continues et la somme de leurs normes est
  strictement inf{\'e}rieure {\`a} 1. Ceci d{\'e}coule
  aussi de la proposition~\ref{continuitedutruc}~: $B~:X_{T} \times
  Y_{T} \rightarrow Y_{T}$ est bicontinue, donc $B(w,\cdot)$ et
  $B(\cdot,w)$ sont dans
  $\mathcal{L}(Y_{T})$, et $\|B(w,\cdot)\|_{\mathcal{L}(Y_T)} +
  \|B(\cdot,w)\|_{\mathcal{L}(Y_T)} < 1$ si $\| w
  \|_{X_{T}}$ est assez petite.
\end{itemize}

Nous disposons maintenant de $v$ solution de $(NS)$ sur $[0,T]$,
appartenant {\`a} $X_{T}\cap Y_{T}$. $\blacksquare$

\bigskip

En particulier, il existe $\tau > 0$ tel que $v(\tau) \in
L^2$ (Parler de $v(\tau)$ a un sens, en effet, nous verrons dans la
suite que $v$ est faiblement continue). Il
nous faut maintenant {\'e}tendre {\`a} $\mathbb{R}^+$ notre solution
d{\'e}finie pour l'instant
localement. C'est l'objet de la proposition suivante, qui porte sur
une estimation d'{\'e}nergie a priori.

\begin{prop}
\label{apriopriori}
Soit $v_{\tau}\in L^2$, et $w$ comme ci-dessus. On consid{\`e}re $v$ solution de
\begin{equation}
\label{troiseq}
 \left\{ \begin{array}{l} 
\partial_{t}v - \Delta v + \mathbb{P}(v \cdot \nabla v +  w \cdot
\nabla v  + v \cdot \nabla w) = 0 \\
\operatorname{div} v = 0 \\
v_{|t=\tau} = v_{\tau}
\end{array} \right.
\end{equation}
On dispose alors de l'estimation a priori suivante sur la norme de $v$
\begin{equation}
\label{inegaal}
\| v(t) \|_{L^{\infty}([\tau,t],L^{2}) ~ \cap
~  L^{2}([\tau,t],\dot{H}^{1})} \leq C \left( \frac{t}{\tau}
\right)^{C \epsilon} \| v_{\tau} \|_{L^2} \, \, .
\end{equation}
\end{prop}

\textsc{Preuve de la proposition~:}

Dans la premi{\`e}re des trois {\'e}quations de~(\ref{troiseq}) on prend le produit
scalaire (spatial) avec~$v$ puis on int{\`e}gre (en
temps) et on obtient, en utilisant $\operatorname{div}
v=\operatorname{div} w=0$~:
\begin{equation}
\label{eqenergie}
\| v(t) \|_{2}^{2} + 2\int_{\tau}^{t} \| \nabla v(s)
\|_{2}^{2}ds + \int_{\tau}^{t} \int_{\mathbb{R}^{2}}(v \cdot
\nabla) v w\, dx\, ds = \| v_{\tau} \|_{2}^{2} \, .
\end{equation}
(Si $p \in [1,\infty]$, on note $\|\cdot \|_p$ pour la norme de
l'espace de Lebesgue $L^p$). L'in{\'e}galit{\'e} de H{\"o}lder et le fait que $\sqrt{t}\| w(t)
\|_{\infty} \leq C \epsilon$ permettent d'{\'e}crire~:
\begin{equation*}
\begin{split}
\left| \int_{\tau}^{t} \int_{\mathbb{R}^{2}}(v \cdot \nabla) v w \, dx
  \, ds
\right| & \leq
\int_{\tau}^{t} \| v(s) \|_{2} \| \nabla v(s)
\|_{2} \| w(s) \|_{\infty} ds \\
& \leq C \epsilon \left( \int_{\tau}^{t} \| \nabla v(s)
  \|_{2}^{2}ds + \int_{\tau}^{t} \frac{\| v(s)
    \|_{2}^{2}}{s}ds \right)
\end{split}
\end{equation*}
En reportant dans~(\ref{eqenergie}), on obtient~:
$$
\| v(t) \|_{2}^{2} + (2 - C\epsilon)\int_{\tau}^{t}
\| \nabla v \|_{2}^{2}~~~~ \leq~~~~ \| v_{\tau}
\|_{2}^{2} + C \epsilon \int_{\tau}^{t} \frac{\| v(s)
  \|_{2}^{2}}{s}ds \, ,
$$
et le lemme de Gronwall permet d'obtenir le r{\'e}sultat souhait{\'e}. $\blacksquare$

\bigskip

On peut maintenant appliquer le sch{\'e}ma standard~: r{\'e}gularisation de
l'{\'e}quation $(NS)$ sur~$[\tau, + \infty[$, puis passage {\`a} la limite faible
en utilisant l'estimation pr{\'e}c{\'e}dente. On en d{\'e}duit que $v$ se
prolonge {\`a} $\mathbb{R}^{+}$ en une solution de $(NS)$
v{\'e}rifiant~(\ref{inegaal}). En posant $u = v+w$, on obtient une
solution globale de $(NS)$ issue de $u_{0}$.
Comme $L^{2}\hookrightarrow \partial BMO$, et du fait des
in{\'e}galit{\'e}s~(\ref{inegaal}) et~(\ref{BdBMO}), on dispose de plus de l'estimation
suivante~:
$$
\| v(t) \|_{\partial BMO} \leq C(\delta) (1+t^\delta) \, \, ,
$$
pour tout $\delta > 0$, o{\`u} l'on a {\'e}crit $C(\delta)$ pour une
constante d{\'e}pendant de $\delta$.

\bigskip

Il nous reste {\`a} d{\'e}montrer la proposition \ref{continuitedutruc}, que nous avons
utilis{\'e}e plus haut~:

\medskip

\textsc{Preuve de la proposition~\ref{continuitedutruc}:}
Les in{\'e}galit{\'e}s de Young et H{\"o}lder permettent d'{\'e}crire
\begin{equation*}
\begin{split}
\| B(w,v)(t) \|_{2} & = \| \int_{0}^{t}
\frac{1}{(t-s)^{3/2}} G \left( \frac{\cdot}{\sqrt{t-s}} \right) \ast
(w(s) v(s)) ds \|_{2} \\
& \leq \int_{0}^{t} \frac{1}{(t-s)^{3/2}} \| G \left(
  \frac{\cdot}{\sqrt{t-s}} \right) \|_{1} \| w(s)
\|_{\infty} \| v(s) \|_{2} ds \\
& \leq \| G \|_{1} \left( \sup_{t>0} \sqrt{t} \| w(t)
\|_{\infty} \right) \left( \sup_{t>0} \| v(t)
\|_{2} \right)
\int_{0}^{t} \frac{ds}{\sqrt{t-s}\sqrt{s}} \, . \\
\end{split}
\end{equation*}
Ainsi,
\begin{equation*}
\begin{split}
\| B(w,v) (t) & \|_{L^{\infty}([0,T],L^{2})} = \sup_{0<t<T}
  \| B(w,v)(t) \|_{2} \\
& \leq \| G \|_{1} \left( \sup_{0<t<T} \int_{0}^{t}
  \frac{ds}{\sqrt{s}\sqrt{t-s}} \right) \left( \sup_{0<t<T} \sqrt{t}
  \| w(t) \|_{\infty} \right) \left( \sup_{0<t<T}
  \| v(t) \|_{2} \right) \\
& \lesssim \| w  \|_{X_{T}} \| v \|_{Y_{T}} \, \, ,
\end{split}
\end{equation*}
o{\`u} l'on note $a \lesssim b$ si $a \leq Cb$ pour une constante $C$.
Il nous reste maintenant {\`a} obtenir une estimation de la norme du
gradient de $B$~; pour ce faire, nous d{\'e}composons cette fonctionnelle
en une somme de deux termes~:
\begin{equation*}
\begin{split}
B(w,v)(t) & = \int_{0}^{t/2}e^{(t-s)\Delta} \mathbb{P} \nabla
(w(s) \otimes v(s))ds + \int_{t/2}^{t}e^{(t-s)\Delta} \mathbb{P}
\nabla (w(s) \otimes v(s))ds \\
& = B_{1}(w,v)(t) + B_{2}(w,v)(t)
\end{split}
\end{equation*}

Et nous {\'e}crivons $\nabla B(w,v) = \nabla B_{1}(w,v) + B_{2}(\nabla
w,v) + B_{2}(w,\nabla v)$. Il vient alors, en utilisant les
in{\'e}galit{\'e}s de Young et H{\"o}lder~:

\begin{equation*}
\begin{split}
\| \nabla B_{1}(w,v)(t) \|_{2} & = \| \int_{0}^{t/2}
\frac{1}{(t-s)^{2}}G(\frac{\cdot}{\sqrt{t-s}})\ast(w(s)v(s))ds
  \|_{2} \\
& \lesssim \int_{0}^{t/2} \frac{ds}{(t-s)\sqrt{s}} \left( \sup_{0<t<T}
    \sqrt{t} \| w(t) \|_{\infty} \right) \left(
    \sup_{0<t<T} \| v(t) \|_{2}  \right) \\
\end{split}
\end{equation*}

Soit~:
$$
\underset{0<t<T}{\sup} \sqrt{t} \| \nabla B_{1}(w,v)(t)
\|_{2} \lesssim
\| w \|_{X_{T}} \| v \|_{Y_{T}} \, \, .
$$
En utilisant encore une fois la m{\^e}me majoration, on montre enfin~:
\begin{equation*}
\begin{split}
\sup_{0<t<T} \| B_{2}(\nabla w,v)(t) \|_{2} & \lesssim
\sup_{0<t<T} \int_{t/2}^{t} \frac{ds}{s \sqrt{t-s}} \left( \sup_{0<t<T}
    t \| \nabla w(t) \|_{\infty} \right) \left(
    \sup_{0<t<T} \| v(t) \|_{2}  \right) \\
& \lesssim \| w \|_{X_{T}} \| v \|_{Y_{T}} \\
\sup_{0<t<T} \| B_{2}(w,\nabla v)(t) \|_{2} & \lesssim
\sup_{0<t<T} \int_{t/2}^{t} \frac{ds}{s \sqrt{t-s}} \left( \sup_{0<t<T}
    \sqrt{t} \| w(t) \|_{\infty} \right) \left(
    \sup_{0<t<T} \sqrt{t} \| \nabla v(t) \|_{2}  \right) \\
& \lesssim \| w \|_{X_{T}} \| v \|_{Y_{T}} \, \, .
\end{split}
\end{equation*}
Il appara{\^\i}{}t maintenant pourquoi nous avons d{\^u} d{\'e}composer $B$
sous la forme $B=B_{1}+B_{2}$. les formules ci-dessus comprennent en
effet des int{\'e}grales (comme $\int_{t/2}^{t} \frac{ds}{s\sqrt{t-s}}$)
qui ne convergeraient pas si le domaine d'int{\'e}gration {\'e}tait
$[0,t]$. $\blacksquare$

\subsection{Preuve de l'unicit{\'e} dans $\mathcal{E}$}

\label{preuvedutheo2}

Cette section est consacr{\'e}e {\`a} la preuve de la partie "unicit{\'e}"
du th{\'e}or{\`e}me~\ref{theoprincipal}.

\medskip

Montrons d'abord que deux solutions dans $\mathcal{E}$ de m{\^e}me
condition initiale sont {\'e}gales. Soient donc $u$, $\tilde{u}$ deux
solutions de $(NS)$ issues de $u_0$ et appartenant {\`a}
$\mathcal{E}$. Soit 
$$\tau = \inf\{t\in \mathbb{R}^+ \, , \, u(t)\neq
\tilde{u}(t)\} \; .$$
Supposons par l'absurde $\tau <\infty$. Il d{\'e}coule alors du lemme
suivant  que
$u(\tau) = \tilde{u}(\tau)$.

\begin{lemm}
Soit $u$ une solution de $(NS)$ telle que $u_{0} \in \partial BMO$ et
$\| u \|_{\mathcal{C},T}$
soit fini. Alors~$u$ est faiblement continue sur $[0,T]$, c'est {\`a}
dire que si $\phi \in \mathcal{S}$, $\operatorname{div} \phi = 0$
$$
<f(t),\phi> \underset{t\rightarrow \tau}{\longrightarrow}
<f(\tau),\phi> \, \, .
$$
\end{lemm}

\textsc{Preuve~:} Soit $\phi \in \mathcal{S}$ et $u$ comme dans
l'{\'e}nonc{\'e}. Puisque $u$ est solution de $(NS)$, elle s'{\'e}crit
$$
u = e ^{t\Delta} u_0 - B(u,u) \, \, .
$$
La continuit{\'e} faible de $e ^{t\Delta} u_0$ est claire. Examinons
maintenant $B(u,u)$. 
\begin{equation}
\label{prodscalBp}
<B(u,u),\phi> = < \int_0^t e^{(t-s)\Delta} \mathbb{P} \nabla u(s)^2 ds , \phi > =
\int_{0}^{t} < u(s)^2 , e^{(t-s)\Delta}\nabla \phi > ds \, \, .
\end{equation}
Comme $\| u \|_{\mathcal{C},T} < \infty$, $u$ v{\'e}rifie,
pour une constante $C$ ind{\'e}pendante de $x \in \mathbb{R}^n$,
\begin{equation}
\label{majorationu} 
\int_0^T {\int\!\!\!\!\!\!-}_{B(x,\sqrt{T})} u(x,t)^2 dx dt \leq C \, \, .
\end{equation}
Comme $\phi$ a une tr{\`e}s forte d{\'e}croissance spatiale, la derni{\`e}re
int{\'e}grale de~(\ref{prodscalBp}) converge bien. Assurons nous
maintenant de la continuit{\'e} faible de $B(u,u)$, en consid{\'e}rant
$$B(u,u)(t) - B(u,u)(t') \, \, ,$$avec $t'<t$.
\begin{equation*}
\begin{split}
<B(u,u)&(t),\phi> - <B(u,u)(t'),\phi> = \\
& \int_{t'}^t <u(s)^2,e
^{(t-s)\Delta} \nabla \phi> ds + \int_0^{t'} <u(s)^2, (e ^{(t-t')\Delta}
- I)e ^{(t'-s)\Delta} \phi > ds
\end{split}
\end{equation*}
En utilisant comme plus haut que $u$ v{\'e}rifie la majoration~(\ref{majorationu}) et
que $\phi \in \mathcal{S}$, on voit que le deuxi{\`e}me membre de cette
derni{\`e}re {\'e}galit{\'e} tend vers 0 si $t'\rightarrow t$, ce qui prouve le lemme. $\blacksquare$

\medskip

Ainsi $u(\tau) = \tilde{u}(\tau)$~; on en d{\'e}duit en
utilisant la proposition~\ref{propunicite} que $u=\tilde{u}$ au
voisinage de~$\tau$~; on contredit ainsi la d{\'e}finition de~$\tau$.

\bigskip

Il nous faut maintenant montrer que les solutions que nous
construisons sont effectivement dans $\mathcal{E}$. Nous consid{\'e}rons
$u=w+v$ une solution construite comme dans le paragraphe
pr{\'e}c{\'e}dent.

Les fonctions $w$ et $v$ v{\'e}rifient chacune la condition $(i)$ de la d{\'e}finition
de $\mathcal{E}$, donc leur somme $u$ aussi.

Reste la condition $(ii)$. Elle est clairement v{\'e}rifi{\'e}e par $w$,
puisque sa
norme dans $X_{\infty}$ est major{\'e}e par $C\epsilon$, et on prend le
param{\`e}tre $\epsilon$ assez petit. Pour montrer que $v$ satisfait aussi
cette condition, nous nous pla{\c c}ons en $\tau > 0$ (le cas $\tau = 0$
est clair) et nous allons voir
que pour un certain $\delta>0$, $v(\tau + \cdot) \in
\mathcal{E}_{\delta}$. Pour plus de simplicit{\'e} dans les notations, nous {\'e}crirons dans la
suite $\tau = 0$. L'{\'e}quation suivante est satisfaite par $v$~:
\begin{equation}
\label{derniereq}
v = e^{t\Delta}v_0 - B(v,v) - B(w,v) -B(v,w) \, \,.
\end{equation}
(noter que $v$ est fortement continue {\`a} valeurs dans $L^2$, donc $v(\tau)=v_0$ est bien d{\'e}finie).
En r{\'e}solvant cette {\'e}quation par point fixe dans $X_T$, on obtient une
solution $v^\star \in X_T$ pour~$T$ assez petit et telle que
$$
\| v^\star \|_{X_T} \lesssim \|e^{t\Delta}v_0\|_{X_T} \underset{T\rightarrow
  0}{\longrightarrow} 0
$$
car $v_0 \in L^2 \subset \overline{\mathcal{S}}^{\partial BMO}$. Il
suffit maintenant de voir que $v^\star = v$ au voisinage de $0$. Ce
sera possible gr{\^a}ce {\`a} un argument de propagation de r{\'e}gularit{\'e} qui
d{\'e}coulera du
\begin{lemm}
Pour $T>0$ les applications suivantes sont continues :
\begin{equation*}
\begin{split}
& B : L^4([0,T],L^4) \times L^4([0,T],L^4) \longrightarrow
L^4([0,T],L^4) \\
& B : X_T \times L^4([0,T],L^4) \longrightarrow
L^4([0,T],L^4) \; . \\
\end{split}
\end{equation*}
\end{lemm}
\textsc{Preuve :} Il est bien connu (voir par
exemple~\cite{bibfabes}) que $B$ est bicontinu sur $L^p L^q$ d{\`e}s
que $\displaystyle \frac{2}{p} + \frac{2}{q} = 1$, ce qui est le cas
ici. Il reste donc seulement {\`a} prouver la seconde assertion du
lemme.
\begin{equation*}
\begin{split}
\|B(u,v)\|_4 & \lesssim \int_0^t \frac{1}{\sqrt{t-s}} \|G\|_1 \|u\|_\infty
\|v\|_4 \, ds \\
& \lesssim \int_0^t \frac{1}{\sqrt{s}\sqrt{t-s}} \|v\|_4 ds
\|u\|_{X_T} \, \, .
\end{split}
\end{equation*}
Or si $t \mapsto \|v(t)\|_4 \in L^4$, les lois de produit entre espaces de Lorentz
impliquent que 
$$
\frac{1}{\sqrt{s}} \| v(s) \|_4 \in L^{4/3,4}
$$
puis que $\displaystyle \int_0^t \frac{1}{\sqrt{s}\sqrt{t-s}} \|v\|_4
\, ds \, \in
L^4$. Autrement dit,
\begin{equation*}
\|B(u,v)\|_4 \lesssim \|v\|_{L^4 L^4} \|u\|_{X_T} \, \, .
\end{equation*}
$\blacksquare$

\smallskip

Par propagation de r{\'e}gularit{\'e} (lemme~\ref{propagation}), $v^\star
\in L^4([0,T],L^4)$. D'autre part, $v$ est pour $T$ fini dans
$L^\infty([0,T],L^2) \cap L^2([0,T],\dot{H}^1)$, donc dans
$L^4([0,T],L^4)$. 
On revient maintenant {\`a} l'{\'e}quation~(\ref{derniereq}), que
l'on r{\'e}soud par point fixe dans $L^4([0,T],L^4)$~; pour $T$ assez
petit, on a unicit{\'e} des solutions de norme petite~; on en d{\'e}duit
que $v^\star = v$ au voisinage de $0$, ce qui conclut la preuve de
l'unicit{\'e} dans $\mathcal{E}$.

\bigskip

On peut maintenant montrer que pour $u_0 \in
\overline{\mathcal{S}}^{\partial BMO}$, la solution~$u$ (unique dans
$\mathcal{E}$) de $(NS)$ est {\`a} valeurs dans 
$\overline{\mathcal{S}}^{\partial BMO}$.

Pour ce faire, soient deux suites $(w_0^\delta)$ et $(v_0^\delta)$
telles que
\begin{equation}
\left\{
\begin{array}{l}
\| w_0^\delta \|_{\partial BMO} \leq \delta \\
v_0^\delta \in \mathcal{S} \\
w_0^\delta + v_0^\delta = u_0
\end{array}
\right.
\end{equation}
et on applique le sch{\'e}ma de construction de solutions du paragraphe
pr{\'e}c{\'e}dent {\`a} la d{\'e}composition $u_0 = w_0^\delta + v_0^\delta$
pour obtenir
$$
u = w^\delta + v^\delta
$$
o{\`u} $v$ est localement born{\'e}e dans $L^2$ et $\| w \|_{\partial BMO}
\leq C \delta$. On conclut en laissant $\delta$ tendre vers 0.

\subsection{Ajout d'une force ext{\'e}rieure}

Nous nous int{\'e}ressons maintenant au syst{\`e}me $(NSF)$. Il admet
lui aussi une formulation int{\'e}grale, sous la forme
\begin{equation*}
u(t) = e^{t\Delta}u_{0} - \int_{0}^{t}e^{(t-s)\Delta} \mathbb{P} \nabla
(u(s) \otimes u(s))ds  + \int_{0}^{t}e^{(t-s)\Delta} \mathbb{P} \nabla
V(s)\,ds 
\end{equation*}
Que doit-on imposer sur $V$ pour que
les r{\'e}sultats obtenus en l'absence de force ext{\'e}rieure soient
toujours valides~?
\begin{itemize} 
\item Si l'on se place dans le cadre des solutions
  d'{\'e}nergie finie ($u_0 \in L^2$), il suffit que $V \in L^{2}(\mathbb{R}^+,L^2)$
  pour garantir l'existence d'une solution.
\item Si maintenant l'on consid{\`e}re le cas $u_0 \in \partial BMO$, on
  peut, comme dans~\cite{bibcannoneplanchon}, examiner quelles sont les
  conditions sur $V$ pour que la m{\'e}thode de point fixe s'applique
  encore et donne une solution globale de $(NSF)$ pour une donn{\'e}e
  initiale assez petite. En reprenant les estimations de~\cite{bibdub}
  sur $\|
  B(u,v)\|_{X_{T}}$, il appara{\^\i}{}t que
$$
\left\| \int_{0}^{t}e^{(t-s)\Delta} \mathbb{P} \nabla V(s)\,ds \right\|_{X_\infty}
\leq \| V \|_{Z}
$$
o{\`u} $\| V \|_{Z}$ est donn{\'e} par
$$
\| V \|_{Z} \overset{\mbox{d{\'e}f}}{=} \sup_{x,R} \int_0^{R^2}
{\int\!\!\!\!\!\!-}_{B(x,R)}
|V(y,t)| dy dt + \sup_{t>0} t \|V(t)\|_{\infty} + \sup_{t>0} t^{3/2}
\| \nabla V(t)\|_{\infty} \, \, .
$$
Ainsi si $\| V \|_{Z}$ et $\|u_{0}\|_{\partial BMO}$ sont assez
petits, la m{\'e}thode de point fixe permet d'obtenir une solution
globale dans $X_{\infty}$ de (NSF).
\end{itemize}

Les deux observations qui viennent d'{\^e}tre faites conduisent au
th{\'e}or{\`e}me suivant~:

\begin{theo}
Soit $u_{0} \in \overline{\mathcal{S}}^{\partial BMO}$. Il existe
$\epsilon>0$ tel que, si
$$
V = V_{v} + V_{w} \mbox{    avec    } \| V_{w} \|_{Z} < \epsilon \mbox{
   et    } V_v \in L^{2}(\mathbb{R}^+,L^2) \; ,
$$
le syst{\`e}me (NSF) admette une solution globale {\`a} valeurs dans
$\partial BMO$ issue de $u_{0}$. De plus, pout tout $\delta$
strictement positif, il existe une constante $C(\delta)$, d{\'e}pendant de $u_0$,
$V$ et $\delta$, telle que
$$
\forall t>0 \, \, \, \, \| v(t) \|_{\partial BMO} \leq C(\delta)(1+t^\delta)
\, .
$$
\end{theo}

\section{Espaces de Besov, $\partial BMO$ et
  $\overline{\mathcal{S}}^{\partial BMO}$}

\subsection{Espaces de Besov}

Une fonction $f$ de $\mathbb{R}^{n}$ dans
$\mathbb{R}^{p}$ est dans l'espace de Besov $\dot{B}^{s}_{p,q}$ si et
seulement si~:
$$ \| f \|_{\dot{B}_{p,q}^{s}} = \left( \sum_{j \in
    \mathbb{Z}} (2^{js} \| \Delta_{j} f \|_{L^{p}})^{q}
\right)^{1/q} < \infty \; ,$$
o{\`u} les $\Delta_{j}$ sont les blocs dyadiques associ{\'e}s {\`a} une
    d{\'e}composition de Littlewood-Paley homog{\`e}ne~:
\begin{gather*}
\Phi \in \mathcal{S} \\
\operatorname{Supp}(\hat{\Phi}) \subset \mathcal{C}(0,3/4,8/3) \\
\Delta_{j} = \hat{\Phi}(2^{-j}D) \\
\sum_{j \in \mathbb{Z}} \Delta_{j} = Id \mbox{ dans } \mathcal{S}' \,
\, .
\end{gather*}
On pourra trouver une pr{\'e}sentation exhaustive des espaces de Besov
dans~\cite{bibtriebel}.

\subsection{L'espace $\partial BMO$}

La r{\'e}f{\'e}rence majeure pour $BMO$ est le livre de
Stein~\cite{bibstein}. Suivant en cela Koch et Tataru, nous
adoptons la d{\'e}finition suivante pour $BMO$~: $f$ appartient {\`a}
$BMO$ si et seulement si~:
$$ \sup_{R>0 , x\in \mathbb{R}^2} \int_{0}^{R^{2}}
{\int\!\!\!\!\!\!-}_{B(x,R)} | \nabla e^{t\Delta} f(x) |^{2} dx dt <
\infty .$$
(o{\`u} l'on note $\displaystyle {\int\!\!\!\!\!\!-}_{B(x,R)} =
\frac{1}{|B(x,R)|}\int_{B(x,R)}$). Il est {\`a} noter (\cite{bibstein},
\cite{bibdub}) que cette d{\'e}finition ne co{\"\i}{}ncide pas
enti{\`e}rement avec la d{\'e}finition "classique" de $BMO$, selon
laquelle $f$ appartient {\`a} $BMO$ si et seulement si~:
$$
\sup_{x \in \mathbb{R}^{d}, R<\rho} {\int\!\!\!\!\!\!-}_{B(x,R)} \left| f(y) -
  {\int\!\!\!\!\!\!-}_{B(x,R)} f \right| dy \, < \, +\infty .
$$
Du moins les fonctions qui v{\'e}rifient la condition ci-dessus sont-elles
dans l'espace $BMO$ tel que nous l'avons d{\'e}fini.

\bigskip

On peut maintenant d{\'e}finir $\partial BMO$, l'espace des
fonctions d{\'e}riv{\'e}es de fonctions de $BMO$. Plus pr{\'e}cis{\'e}ment, on
dira que $f \in \partial BMO(\mathbb{R}^{n})$ si et seulement si il existe des
fonctions de~$BMO(\mathbb{R}^{n})$ $f_{1},\dots,f_{n}$ telles que~: 
$$
\label{defderiv}
f = \sum_{i=1}^{n} \frac{\partial}{\partial x_{i}} f_{i} \, .
$$
On peut montrer (\cite{bibkoch}) que cette d{\'e}finition est
{\'e}quivalente au fait que $f$ v{\'e}rifie~:
\begin{equation} \label{cestclair} \sup_{R>0 , x\in \mathbb{R}^2}
 \int _{0}^{R^2} {\int\!\!\!\!\!\!-}_{B(x,R)} | e^{t\Delta} f(x)
  |^{2} dx dt < \infty \, \, . \end{equation}
Cette derni{\`e}re formule d{\'e}finit une norme sur $\partial BMO$~; nous
 la noterons $\| \cdot \|_{\partial BMO}$. Muni de cette norme,
 $\partial BMO$ est un espace complet.

Pour finir, notons le lemme suivant~:
\begin{lemm}
\label{chaleurdBMO}
Il existe $C>0$ tel que, si $f \in \partial BMO$ et $t>0$~:
$$
\| e^{t\Delta} u \|_{\partial BMO} \leq \| u \|_{\partial BMO} \, \, .
$$
\end{lemm}

\subsection{Donn{\'e}es initiales permettant d'appliquer le th{\'e}or{\`e}me}

Le th{\'e}or{\`e}me~(\ref{theoprincipal})
est {\'e}nonc{\'e} pour $u_{0} \in \overline{\mathcal{S}}^{\partial BMO}$, mais la
preuve donn{\'e}e plus haut en dit un peu plus. Elle repose en effet sur
un d{\'e}coupage de la condition initiale en une partie r{\'e}guli{\`e}re et une
partie dans $\partial BMO$ de norme $<\epsilon$. On  en d{\'e}duit facilement
qu'il existe $\epsilon$ tel que la m{\^e}me m{\'e}thode s'applique pour $u_{0}
\in \partial BMO$ tel que~:
$$
d(u_{0},\mathcal{S})=\inf_{f\in
  \mathcal{S}}\| f-u_{0} \|_{\partial BMO} <\epsilon \, .
$$ 
Pour $\epsilon$ assez petit, les conclusions du
th{\'e}or{\`e}me~\ref{theoprincipal} restent donc inchang{\'e}es, {\`a} ceci pr{\`e}s
que la solution construite n'est plus {\`a} valeurs dans $\overline{\mathcal{S}}^{\partial BMO}$.

\subsection{L'espace $\overline{\mathcal{S}}^{\partial BMO}$}

Soit $VMO$ (pour \textit{Vanishing Mean Oscillations}) l'espace d{\'e}fini par~:
\begin{equation}
\label{defVMO}
f \in VMO \Leftrightarrow 
\left\{ \begin{array}{l} f \in BMO \\ \displaystyle \sup_{x \in
      \mathbb{R}^{d}, R<\rho} {\int\!\!\!\!\!\!-}_{B(x,R)} \left| f(y) -
  {\int\!\!\!\!\!\!-}_{B(x,R)} f \right| dy \underset{\rho \rightarrow
  0}{\longrightarrow} 0 \\
\displaystyle \, \sup_{x \in
      \mathbb{R}^{d}, R>\rho} {\int\!\!\!\!\!\!-}_{B(x,R)} \left| f(y) -
  {\int\!\!\!\!\!\!-}_{B(x,R)} f \right| dy \underset{\rho \rightarrow
  +\infty}{\longrightarrow} 0 \end{array} \right.
\end{equation}
On trouve dans~\cite{bibstein} le r{\'e}sultat suivant~: 

\begin{prop}
\label{adherenceS}
$VMO$ est l'adh{\'e}rence de $\mathcal{S}$ dans $BMO$.
\end{prop}

Nous notons maintenant $\partial_i$ pour
$\frac{\partial}{\partial x_i}$ et $\partial X$ pour l'ensemble des d{\'e}riv{\'e}es
de fonctions de l'espace $X$ (au sens de~(\ref{defderiv})). La
proposition suivante va nous permettre de mieux appr{\'e}hender l'espace
$\overline{\mathcal{S}}^{\partial BMO}$~: 

\begin{prop}
On a~:
$$ \partial VMO \subset \overline{\mathcal{S}}^{\partial BMO} $$
\end{prop}
\textsc{Preuve de la proposition~:}
Soit $f \in VMO = \overline{\mathcal{S}}^{BMO}$, montrons que
$\partial_{i} f \in 
\overline{\mathcal{S}}^{\partial BMO}$. Soit donc $(\phi_{n})$ une suite
de $\mathcal{S}$ telle que
$$
\phi_{n} \overset{BMO}{\longrightarrow} f
$$
Alors, comme $\| \cdot \|_{BMO} \geq \|\partial_i \cdot \|_{\partial
  BMO}$, on a aussi~:
$$
\partial_i \phi_{n} \overset{\partial BMO}{\longrightarrow} \partial_i f
$$
On en d{\'e}duit donc que~:
$$
\partial_{i} \overline{\mathcal{S}}^{BMO} \subset
\overline{\mathcal{S}}^{\partial BMO},
$$
d'o{\`u} le r{\'e}sultat.
$\blacksquare$

\medskip

Soulignons enfin le ph{\'e}nom{\`e}ne suivant: nous venons de voir que les
fonctions de $BMO$ qui sont dans l'adh{\'e}rence de $\mathcal{S}$ sont
plus r{\'e}guli{\`e}res et ont une meilleure d{\'e}croissance {\`a} l'infini
que les autres fonctions de $BMO$. Ceci est exprim{\'e} par les
conditions dans~(\ref{defVMO}) lorsque $\rho \rightarrow 0$ ou $\rho
\rightarrow \infty$.

Pour $\overline{\mathcal{S}}^{\partial BMO}$, on observe aussi ce
double effet de r{\'e}gularisation locale et de d{\'e}croissance plus
forte {\`a} l'infini. Ainsi, on
v{\'e}rifiera facilement, en une dimension, que $x \mapsto \sin(x)$ et
$x \mapsto \frac{1}{|x|}$ appartiennent {\`a}
$\partial BMO \setminus \overline{\mathcal{S}}^{\partial BMO}$, la premi{\`e}re fonction ne
pr{\'e}sentant pas de d{\'e}croissance {\`a} l'infini et la seconde ayant
une trop forte discontinuit{\'e} en $0$.

\section{Conclusion}

Pour conclure, il est int{\'e}ressant de r{\'e}capituler les r{\'e}sultats
dont nous disposons pour l'{\'e}quation (NS), suivant $u_0$~:
\begin{itemize}
\item Si $u_0$ est dans $\overline{\mathcal{S}}^{\partial BMO}$, il
  existe une solution $u$, {\`a} valeurs dans
  $\overline{\mathcal{S}}^{\partial BMO}$ telle que
$$
\| u(t) \|_{\partial BMO} \leq C(u_0)(1+\sqrt{t}) \; .
$$
\item Si $u_0$ est dans $\partial BMO$ avec $\| u_0 \|_{\partial
  BMO}<\epsilon$, il existe une solution $u$,
  {\`a} valeurs dans~$\partial BMO$, et telle que
$$
\| u(t) \|_{\partial BMO} \leq C \|u_0 \|_{\partial BMO} \, \, .
$$
\item Dans le cas g{\'e}n{\'e}ral $u_0 \in \partial BMO$, on ne sait rien
  sur l'existence d'une solution.
\end{itemize}

\section{Annexe~: r{\'e}sultats de type point fixe}

Nous utilisons dans cet article un th{\'e}or{\`e}me de point fixe et un
lemme de propagation de r{\'e}gularit{\'e}. Ils sont rappel{\'e}s ici (on pourra
se reporter {\`a}~\cite{bibgallag2})~:

\begin{theo}[Existence et unicit{\'e}]
\label{pointfixe}
Soit $X$ un espace de Banach, $L$ un op{\'e}rateur lin{\'e}aire sur $X$ de norme
$\lambda <1$, et $B$ un op{\'e}rateur bilin{\'e}aire tel que~:
$$
\| B(x,y) \|_{X} \leq \gamma \| x \|_{X}
\| y \|_{X}
$$
Alors pour tout $y \in X$ tel que
$$
4\gamma \| y \|_{X} < (1-\lambda)^{2}
$$
la suite d{\'e}finie par~:
$$
\left\{
\begin{array}{ll} X_{0} = 0 \\ X_{n+1} = y + LX_{n} + B(X_{n},X_{n})
\end{array}
\right.
$$
converge dans $X$ vers l'unique solution de
$$
x = y + Lx+B(x,x)
$$
telle que
$$
2\gamma\| x \|_{X}<(1-\lambda)
$$
De plus, on a l'estimation suivante~:
$$
\| x \|_{X} \lesssim \| y \|_{X}
$$
\end{theo}

\begin{lemm}[Propagation de r{\'e}gularit{\'e}]
\label{propagation}
Soient $x$, $(X_{n})_{n\in \mathbb{N}}$, $\lambda$ et $\gamma$ comme dans le th{\'e}or{\`e}me
pr{\'e}c{\'e}dent. On suppose de plus que $y$ est dans un espace de Banach
$Y$, que $L$ est un op{\'e}rateur lin{\'e}aire sur $Y$ de norme $\mu$ et que
\begin{equation*}
\begin{split}
& \| B(f,g) \|_{Y} < \kappa \| f \|_{X}
\| g \|_{Y} \\
& \| B(g,f) \|_{Y} < \kappa \| f \|_{X}
\| g \|_{Y} . \\ 
\end{split}
\end{equation*}
On suppose enfin que $\mu<1$ et que
$\kappa(1-\lambda)<(1-\mu)\gamma$. Alors $X_{n}$ converge vers $x$ dans $Y$ et
$$
\| x \|_{Y} \lesssim \| y \|_{Y}
$$
\end{lemm}

\end{document}